\documentclass{amsart}

\newtheorem{theorem}{Theorem}[section]
\newtheorem{lemma}[theorem]{Lemma}
\newtheorem{cor}[theorem]{Corollary}

\theoremstyle{definition}

\theoremstyle{remark}

\numberwithin{equation}{section}

\font\bbb=msbm10 scaled 1100
\newcommand{\eps}{\epsilon}

\newcommand{\del}{\partial}

\newcommand{\ra}{\rightarrow}

\newcommand{\real}{\mbox{\bbb R}}

\newcommand{\inv}{^{-1}}

\newcommand{\grad}{\nabla}

\renewcommand{\div}{\nabla\cdot}

\newcommand{\Lie}{{\mathcal L}}

\newcommand{\rest}[2]{\left. #1\right\vert_{#2}}		
\newcommand{\be}{\begin{equation}}
\newcommand{\ee}{\end{equation}}
\newcommand{\bea}{\begin{eqnarray}}
\newcommand{\eea}{\end{eqnarray}}
\newcommand{\bmini}{\footnotesize\begin{center}\begin{minipage}{5.5in}}
\newcommand{\emini}{\end{minipage}\end{center}\normalsize}
\newcommand{\pf}{{\em Proof: }}

\newcommand{\eg}{{\em e.g.}}
\newcommand{\ie}{{\em i.e.}}

\begin{document}

\title{STRATIFIED INTEGRALS AND UNKNOTS IN INVISCID FLOWS}

\author{John B. Etnyre}
\address{Department of Mathematics, Stanford University, 
	Palo Alto, CA 94305 }
\email{etnyre@math.stanford.edu}
\thanks{JBE supported in part by NSF Grant DMS-9705949.}

\author{Robert W. Ghrist}
\address{School of Mathematics, Georgia Institute of Technology, Atlanta,
	GA 30332-0160}
\email{ghrist@math.gatech.edu}



\begin{abstract}
We prove that any steady solution to the $C^\omega$ Euler
equations on a Riemannian $S^3$ must possess a periodic
orbit bounding an embedded disc. 
One key ingredient is an extension of Fomenko's work on the
topology of integrable Hamiltonian systems to a degenerate 
case involving stratified integrals. The result on the Euler 
equations follows from this when 
combined with some contact-topological perspectives and a 
recent result of Hofer, Wyzsocki, and Zehnder.
\end{abstract}

\maketitle

\section{Introduction}

The mathematical approach to knot theory initiated by Lord Kelvin
began as a problem in fluid dynamics: to understand the manner in 
which closed flowlines in the \ae ther are partitioned into various
knot types, with the goal of recovering the periodic table \cite{Tho69}. 
Unfortunately, the two subjects quickly diverged and have not 
since come into such close companionship. There are several 
key exceptions of which we mention two. The continual work of 
Moffatt \cite{Mof85,Mof86,Mof94} from the 
engineering side to recognize the role 
that topology plays in physical fluid dynamics has been largely 
responsible for the acceptance of the definition of 
{\em helicity} (an important topological invariant) within the 
applied community. As an example of reversing this scenario, we
mention the work of Freedman and He \cite{FH91a,FH91b}, who
use a physical notion (hydrodynamical energy) 
to define a topological invariant for knots and links. 

There are several unavoidable problems in the attempt to reconcile
knot theory and fluid dynamics, not the least of which is that the 
fundamental starting point, the global existence of solutions to the 
Euler and Navier-Stokes equations on $\real^3$, is unknown and
perhaps not true. Coupled with this difficulty is the fact that viscosity, 
unusual boundary conditions, and poorly-understood phenomenon
of turbulence conspire to make it
nearly impossible to rigorously analyze the solutions to the relevant 
equations of motion, even with powerful analytical techniques
currently available. 

However, since so little is known about the rigorous behavior of 
fluid flows, any methods which can be brought to bear to prove
theorems about their behavior are of interest and of potential 
use in further understanding these difficult problems.
We propose a view of the relevant equations of motion for inviscid
(without viscosity) fluid flows which sets up the possibility of a 
topological approach. We do so by not only restricting the class of 
flows considered (steady, nonsingular) but also by expanding the 
class via ``forgetting'' all information about the metric structure.

\subsection{The Euler equations}

For a recent topological treatment of the 
equations of motion for a fluid, see the excellent 
monograph of Arnold and Khesin \cite{AK98}.
Any discussion of fluid dynamics must begin with the relevant 
equations of motion, the most general of which is the Navier-Stokes 
equation. Let $u$ denote a time-dependent vector field on 
$\real^3$ (the {\em velocity} of the fluid), $p$ denote a time-dependent
real-valued function on $\real^3$ (the {\em pressure}), and $\nu\geq 0$
denote a constant (the {\em viscosity}). Then, the Navier-Stokes 
equations are 
\begin{equation}
	\frac{\del u}{\del t} + (u\cdot\nabla)u 
	- \nu\Delta u = -\nabla p + f \,\,\,\, ; \,\,\,\, \div{u} = 0, 
\end{equation}
where $f$ is a time-dependent function on $\real^3$ representing 
body forces such as gravity, etc. 
The Euler equations (the form that we will concern ourselves with) 
are precisely the Navier-Stokes equations in the 
absence of viscosity and body forces:
\begin{equation}
\label{eq_preEuler}
\frac{\del u}{\del t}+(u\cdot\nabla)u = -\nabla p 
\,\,\,\,\,\, ; \,\,\,\,\,\, \div{u}=0 .
\end{equation}
We will for the remainder of this paper work under the following 
three simplifying assumptions:
\begin{enumerate}
\item
	All flows considered will be inviscid: $\nu=0$.
\item
	All flows considered will be steady: $\frac{\del u}{\del t}=0$.
\item
	All flows considered will be nonsingular: $u\neq 0$ anywhere.
\end{enumerate}
The first step in embedding this problem into a topological setting is
to expand the class of fluid flows that we consider. First, instead of
restricting the flows to $\real^3$, we will allow for fluid flows on 
any three-manifold $M$. In order to make sense of the various 
operations (grad, curl) in Equation~(\ref{eq_preEuler}), we must
choose a Riemannian metric, $g$, with respect to which these 
operations are taken. Finally, in order to work with the volume-preserving
condition of the Euler equations, we must choose an appropriate
volume form, $\mu$. One can of course choose the precise volume 
form $\mu_g$ induced by the metric; however, for 
the sake of generality, we allow for arbitrary $\mu$. This has 
physical significance, as noted in \cite{AK92}. 

The form which the Euler equations now take is the following:
\begin{equation}
\frac{\del u}{\del t} + \nabla_uu = -\grad p 
\,\,\,\,\,\, ; \,\,\,\,\,\, \Lie_u\mu=0,
\end{equation}
where $\nabla_u$ is the covariant derivative along $u$ defined by the 
metric $g$, and $\Lie_u$ is the Lie derivative along $u$. By a suitable 
identity for the covariant derivative (see \cite[p. 588]{AMR88})
one can transform the previous equation into an exterior differential
system:
\begin{equation}\label{eq_Euler}
\frac{\displaystyle \del(\iota_ug)}{\displaystyle \del t} + 
	\iota_ud\iota_ug = -dP \,\,\,\,\, ; \,\,\,\,\, 
	\Lie_u\mu=0 .
\end{equation}
Here, $\iota_ug$ denotes the one-form obtained from $u$ via contraction 
into the first slot of the metric, and $P$ is a modified 
pressure function from $M$ to the reals. 
It is this form of the Euler equation with which we will be 
concerned for the remainder of this paper.

\subsection{Statement of results}

In an earlier paper \cite{EG:beltrami}, the authors initiated the use of 
contact-topological ideas in the study of the Euler equations.
There, it was shown that all steady solutions of sufficiently
high regularity on $S^3$ possess a closed flowline: \ie, the
Seifert Conjecture is true in the hydrodynamical context. This
opens the possibility for asking questions about knotting and
linking phenomena common to all fluid flows on $S^3$.

The main result we prove in this note is the following:
\begin{theorem}
\label{thm_Main}
Any steady solution to the $C^\omega$ Euler equations on a 
Riemannian $S^3$ must possess a closed flowline which bounds
an embedded disc.
\end{theorem}

The proof of this theorem relies upon deep results due to 
several authors, most especially the work of Hofer et al. 
\cite{HWZ96b} on unknotted orbits in Reeb fields, as well as 
the theorem of Wada \cite{Wad89} on nonsingular Morse-Smale
flows. A key ingredient is a generalization in \S\ref{sec_Unknots}
of a theorem of Fomenko and Nguyen \cite{FN91} and 
of Casasayas et al \cite{CMAN} to a degenerate case:

\begin{theorem}
Any nonsingular vector field on $S^3$ having a $C^\omega$ integral of 
motion must possess a pair of unknotted closed orbits. 
\end{theorem}

This theorem has certain peripheral implications in the Fomenko-style
approach to two degree-of-freedom integrable Hamiltonian systems. 
We elaborate upon these themes in \S\ref{sec_Fomenko}.

We note that Theorem~\ref{thm_Main} is but one 
small piece of data concerning knot theory within hydrodynamics
(note in particular the work of Moffatt et al. \cite{Mof94}). 
In a future paper \cite{EG:eulerknot}, we will consider 
the other end of the spectrum: namely, what is possible as 
opposed to what is inevitable. There, we will construct steady
nonsingular solutions to the Euler equation possessing knotted 
orbits of all possible knot types simultaneously. 

\section{Contact / integrable structures for steady Euler flows}

After providing a brief background on contact geometry/topology, 
we consider the class of steady Euler flows on the three-sphere
in the $C^\omega$ category. We demonstrate that there is a 
dichotomy between integrable solutions, and solutions which are
related to contact forms. 

\subsection{Contact structures on three-manifolds}

A more thorough introduction to the field of contact topology
can be found in the texts \cite{Aeb94,MS95}. A {\sc contact structure}
on an odd-dimensional manifold is a completely nonintegrable 
hyperplane distribution. We will restrict to the case of 
a three-manifold. In this case, a contact structure is a completely
nonintegrable  smoothly
varying field of 2-dimensional subspaces of the tangent spaces.
Unlike vector fields, such plane fields do not necessarily 
integrate to form a foliation. This integrability of a plane 
field $\xi$ is measured by
the Frobenius condition on a defining (local) 1-form $\alpha$.
If $\xi = \ker\alpha$, then $\xi$ is a contact structure 
if and only if $\alpha\wedge d\alpha$ vanishes nowhere. Such 
a form $\alpha$ is a {\sc contact form} for $\xi$. 

The topology of the structures $\xi$ and the geometry of the 
associated forms $\alpha$ has of late been a highly active and 
exciting field. As contact structures are in a strong sense the 
odd-dimensional analogue of a symplectic structure \cite{MS95,Arn80}, 
many of the interesting phenomena of that discipline carry over. 
Of particular importance is the existence of a certain class of 
contact structures -- the {\sc tight} structures -- which possess
topological restrictions not otherwise present. Such structures 
are fairly mysterious: basic questions concerning existence and 
uniqueness of such structures are as yet unanswered. See 
\cite{Eli89,ET97} for more information.

Contact geometry has found applications in a number of disciplines.
We note in particular the utility of contact geometry in 
executing a form of reduction in dynamics with special symmetry
properties \cite{HM97}. In a recent work \cite{EG:beltrami}, 
the authors initiated the use of modern contact-topological 
methods in hydrodynamics. Most of these dynamical applications 
revolve around the notion of a {\sc reeb field} for a contact
form. The Reeb field associated to a contact form $\alpha$ is 
the unique vector field $X$ satisfying the equations:
\begin{equation}
	\iota_X\alpha = 1 \,\,\,\,\, ; \,\,\,\,\, 
	\iota_Xd\alpha = 0 .
\end{equation}
The Reeb field forms a canonical section of the {\sc characteristic
line field} $\ker d\alpha$ on $M$. Reeb fields are by definition 
nonsingular and preserve the volume form $\alpha\wedge d\alpha$. 
In \S\ref{sec_Hofer}, we review a result on the topology of Reeb
fields due to Hofer et al.

\subsection{The dichotomy for steady Euler flows}

The following theorem is a specialized version of the 
general correspondence between solutions to the Euler equation 
and Reeb fields in contact geometry derived in \cite{EG:beltrami}.
We include the [simple] proof for completeness.

\begin{theorem}
\label{thm_Dichotomy}
Let $u$ denote a steady nonsingular solution to the Euler 
equations of class $C^\omega$ on a Riemannian $S^3$. Then at least one 
of the following is true:
\begin{enumerate}
\item
	There exists a nontrivial integral for $u$; or
\item
	$u$ is a nonzero section of the characteristic line field
	of a contact form $\alpha$.
\end{enumerate}
\end{theorem}
\pf
If $u$ is a steady solution then 
\begin{equation}
	\iota_ud\iota_u g = -dP .
\end{equation}
As all the data in the equation is assumed real-analytic,
the differential $dP$ must be $C^\omega$ and hence vanishes identically
if and only if it vanishes on an open subset on $S^3$. Note that
\begin{equation}
	\Lie_uP = \iota_udP = \iota_u\iota_u(d\iota_ug) \equiv 0 ,
\end{equation}
and thus a nonconstant $P$ yields a nontrivial 
integral for the vector field $u$.\footnote{This is a ``modern'' 
reformulation of the classical Bernoulli Theorem for fluids.} 

In the case where $dP\equiv 0$, we have that 
$\iota_ud\iota_ug\equiv 0$. Consider the nondegenerate 1-form 
$\alpha:=\iota_ug$ dual to the vector field $u$ via the metric. 
In addition, denote by $\beta:=\iota_u\mu$ the 2-form obtained 
by pairing $u$ with the volume form $\mu$. Since $u$ is $\mu$-preserving,
it is the case that $d\beta\equiv 0$. Also, by definition, 
$\beta$ has a one-dimensional kernel spanned by $u$: $\iota_u\beta=0$.
As $\iota_ud\iota_ug = 0$, it follows that $d(\iota_ug)=h\iota_u\mu$ 
for some function $h:S^3\ra\real$.

It is a classical fact that $h$ is an integral of $u$: observe that 
\begin{equation}
0 = d^2\alpha = d(h\beta) = dh\wedge\beta 
	+ hd\beta = dh\wedge\beta ,
\end{equation}
which implies that $\iota_udh = 0$.

Thus, the only instance in which $u$ is not integrable is when 
$h$ is constant. If $h\equiv 0$,
then by the Frobenius condition, the 1-form $\iota_ug$ defines 
a $C^\omega$ codimension-one foliation of $S^3$ which is transverse to $u$.
This is impossible due to the presence of a $T^2$ leaf in the 
foliation (guaranteed by Novikov's Theorem \cite{Nov67}) which is
transverse to the volume-preserving flow of $u$. Or, alternatively, 
the $C^\omega$ codimension-one foliation of $S^3$ violates the 
Haefliger Theorem \cite{Hae56}. 

Thus, if $u$ is not integrable via $h$, then $h$ is a nonzero 
constant. So, the 1-form $\iota_ug$ dual to $u$ is a contact form
since $\alpha\wedge d\alpha = h\alpha\wedge\beta\neq 0$. Note that in 
this case $u$ is a section of the characteristic line field of
this contact form: $\iota_u(d\iota_ug)=0$.
\qed

\subsection{The results of Hofer et al.}
\label{sec_Hofer}

The recent deep work of Hofer \cite{Hof93} and of 
Hofer, Wyszocki, and Zehnder \cite{HWZ96,HWZ96b} utilizes 
analytical properties of pseudoholomorphic curves in products 
of a contact manifold with $\real$ to elucidate the dynamics 
and topology of Reeb fields. As such, these results are 
directly applicable to the understanding of steady nonsingular
Euler flows. The specific result we employ for this paper 
is the following:

\begin{theorem}[Hofer et al. \cite{HWZ96b}]
Let $\alpha$ be a contact form on a homology 3-sphere $M$. Then the 
Reeb field associated to $\alpha$ possesses a periodic orbit which 
bounds an embedded disc. 
\end{theorem}

This theorem, combined with Theorem~\ref{thm_Dichotomy}, yields a 
proof of the main result [Theorem~\ref{thm_Main}] in the difficult case 
where the velocity field is a section of a characteristic line 
field for a contact form. The integrable case must yet be 
considered.

\section{Knotted orbits in integrable systems}

Having dispensed with the nonintegrable cases, we turn in this section 
to consider the knot types associated to periodic orbits in 
integrable Hamiltonian systems. Several fundamental results about
such flows are presented in the work by Casasayas et al. 
\cite{CMAN} and by Fomenko and Nguyen \cite{FN91}: we will 
review these results and extend them to the degenerate framework
we require in order to complete the proof of the main theorem.

Recall that a two degree-of-freedom Hamiltonian system on a 
symplectic four-manifold $(W,\omega)$ with Hamiltonian $H$ 
splits into invariant codimension-one submanifolds $Q_c = H\inv(c)$ 
(on the regular values of $H$). The system is said to be {\sc integrable}
on $Q_c$ if there exists a function $F:Q_c\ra\real$ such that 
$F$ is independent of $H$ and $\{F,H\}=0$: in other words, $F$ is invariant
on orbits of the Hamiltonian flow. 

A topological classification of such two degree-of-freedom integrable
systems exists when the integrals are nondegenerate. An integral $F$ is 
said to be {\sc bott} or {\sc bott-morse} if the critical point set $cp(F)$
consists of a finite collection of 
submanifolds $\Sigma_i$ of $Q_c$ which are transversally 
nondegenerate: the restriction of the Hessian $d^2F$ to the 
normal bundle $\nu\Sigma_i$ of $cp(F)$ is nondegenerate. Such 
integrals are generic among integrable systems, but given a 
particular system, it is by no means easy to verify whether an 
integral is Bott-Morse.

The classification of knotted orbits in Bott-integrable 
Hamiltonian systems on $S^3$ is best accomplished via the 
classification of round-handle decompositions of $S^3$. 
Recall that a {\sc round handle} in dimension three is a solid torus
$H=D^2\times S^1$ with a specified index and exit set
$E\subset T^2=\del(D^2\times S^1)$ as follows:
\begin{description}
\item[index 0] $E = \emptyset$.
\item[index 1] $E$ is either (1) a pair of disjoint annuli on the boundary
	torus, each of which wraps once longitudinally; or 
	(2) a single annulus which wraps twice longitudinally.
\item[index 2] $E=T^2$. 
\end{description}
A {\sc round handle decomposition} (or RHD) for a manifold $M$ is a 
finite sequence of submanifolds
\begin{equation} 
	\emptyset = M_0\subset M_1\subset\cdots M_n=M ,
\end{equation}
where $M_{i+1}$ is formed by adjoining a round handle 
to $\del M_i$ along the exit set $E_{i+1}$ of the round handle.
The handles are added in order of increasing index. 
Asimov \cite{Asi75} and Morgan \cite{Mor78} 
used round handles to classify nonsingular
{\sc morse-smale} vector fields: that is, vector fields whose
recurrent sets consist entirely of a finite number of hyperbolic
closed orbits with transversally intersecting invariant
manifolds. In short, the cores of an RHD, labelled by the index,
correspond to the periodic orbits of a nonsingular Morse-Smale
vector field, labelled by the Morse index \cite{Asi75}. 
The classification of RHD's (in the context of nonsingular Morse-Smale 
flows) on $S^3$ was achieved by Wada \cite{Wad89} following
work of Morgan \cite{Mor78}. 
\begin{theorem}[Wada \cite{Wad89}]
\label{thm_Wada}
Let ${\mathcal W}$ be the collection of indexed links determined by the
following eight axioms: 
\begin{itemize}
\item [O] The Hopf link indexed by 0 and 2 is in ${\mathcal W}$.
\item [I] If $L_{1}, L_{2} \in \mathcal W$ then $L_{1} \circ L_{2} \circ u
\in \mathcal W$, where $u$ (here and below) is an unknot in $S^{3}$
indexed by 1, and $\circ$ denotes the {\sc split sum} of knots (\ie, 
separable by means of an embedded 2-sphere).
\item [II] If $L_{1}, L_{2} \in \mathcal W$ and $K_{2}$ is a component of
$L_{2}$ indexed by 0 or 2, then $L_{1} \circ (L_{2} -K_{2}) \circ u
\in \mathcal W$.
\item [III] If $L_{1},L_{2} \in \mathcal W$ and $K_{1},K_{2}$ are components
of $L_{1},L_{2}$ with indices 0 and 2 (resp.), then $(L_{1}-K_{1}) \circ
(L_{2}-K_{2}) \circ u \in \mathcal W$.
\item [IV] If $L_{1},L_{2} \in \mathcal W$ and $K_{1},K_{2}$ are components
of $L_{1},L_{2}$ (resp.) each with index 0 or 2, then 
\[((L_{1},K_{1}){\#}(L_{2},K_{2})) \cup m \in \mathcal W, \]
where $K_1\# K_2$ shares the
index of either $K_1$ or $K_2$ and $m$ is a meridian of $K_{1} \# K_{2}$
indexed by 1.
\item [V] If $L \in \mathcal W$ and $K$ is a component of $L$ indexed 
by $i=0$ or 2, then $L' \in \mathcal W$, where $L'$ is obtained 
from $L$ replacing a tubular neighborhood of $K$ with a solid torus
with three closed orbits, $K_1$, $K_2$, and $K_3$. $K_1$ is the core 
and so has the same knot type as $K$. $K_2$ and $K_3$ are parallel 
$(p,q)$ cables
of $K_1$. The index of $K_2$ is 1. The indices of $K_1$ and $K_3$ may 
be either 0 or 2, but at least one of them must be equal to the index
of $K$.
\item [VI] If $L \in \mathcal W$ and $K$ is a component of $L$ indexed 
by $i=0$ or 2, then $L' \in \mathcal W$, where $L'$ is obtained 
from $L$ by changing the index of $K$ to 1 and placing
a $(2,q)$-cable of $K$ in a tubular neighborhood of $K$,
indexed by $i$.
\item [VII] ${\mathcal W}$ is minimal. That is, ${\mathcal W} \subset {\mathcal 
W}'$
for any collection, ${\mathcal W}'$, satisfying O-VI.
\end{itemize}
Then the class of indexed periodic orbit links arising within 
nonsingular Morse-Smale flows on $S^3$ is precisely ${\mathcal W}$. 
\end{theorem}

\begin{cor}[Wada \cite{Wad89}]
Every smooth nonsingular Morse-Smale vector field on $S^3$ 
possesses a pair of unknotted closed orbits.
\end{cor}
\pf
The base Hopf link is such a pair. It is clear that the Wada moves
I-VI leave this property invariant.
\qed

The relationship between Wada's Theorem and nonsingular integrable 
Hamiltonian systems was developed by Casasayas et al. \cite{CMAN}.
The idea is straightforward: given a Bott-integrable nonsingular Hamiltonian 
system with integral $P$, the vector field $-\grad P$ is a field 
with curves of Bott-Morse type critical points. A small perturbation 
tangent to the critical curves yields a nonsingular Morse-Smale flow.
\begin{cor}[Casasayas et al. \cite{CMAN}]
\label{cor_CMAN}
Every Bott-integrable $C^\infty$ Hamiltonian flow on a symplectic
4-manifold having a nonsingular $S^3$ energy surface possesses a pair of 
unknotted invariant critical curves. 
\end{cor}
Similar results were obtained by
Fomenko and Nguyen \cite{FN91}. 

\section{Proof of Theorem \ref{thm_Main}}
\label{sec_Unknots}

As is clear from the previous section, we may obtain information 
about the knot data of integrable Euler fields if we can ensure
that the critical sets are all of Bott-Morse type. 
Nowhere in the literature is there a discussion of the non-Bott case 
with respect to knotting and linking phenomena. This is very
difficult if not impossible to control in the general $C^\infty$ case; 
however, in the real-analytic case, we may still analyze the 
degenerate critical point sets.

\begin{lemma}
\label{lem_Strat}
Any critical set of a nontrivial $C^\omega$ integral $P$ for a 
nonsingular vector field $X$ on $S^3$ is a [Whitney] 
stratified set of (topological) dimension at most two.
\end{lemma}
\pf
Denote by $cp(P)$ the critical points of $P$ and by $\Sigma$ 
a connected component of the inverse image of the critical values of $P$.
It follows from the standard theorems concerning real-analytic 
varieties \cite{Whi57,GM88} that the set $\Sigma$ is a (Whitney) 
stratified set. That is, although $\Sigma$ is not a manifold, 
it is composed of manifolds --- or {\sc strata} --- glued together
along their boundaries in a controlled manner. It follows from 
analyticity that $\Sigma$ has topological dimension 
less than or equal to two; otherwise, $P$ would be a constant.
\qed

\begin{lemma}
\label{lem_Branched}
The critical set $\Sigma$ is either an embedded closed 
curve in $S^3$, or else is a 
(non-smoothly) branched 2-manifold, where the non-manifold set
of $\Sigma$ is a finite invariant link in $S^3$. The complement of 
this set in $\Sigma$ [the 2-strata] consists of critical tori, 
as well as annuli and M\"obius bands glued to the singular 
link along their boundaries.
\end{lemma}
\pf
As $\rest{X}{\Sigma}$ is a nonsingular
vector field, $\Sigma$ must have a stratification devoid of 
0-strata: only 1- and 2-strata are permitted. Furthermore, 
the topology of $\Sigma$ must be transversally homogeneous
with respect to the flow: a neighborhood of any point $x$ in 
$\Sigma$ is homeomorphic to a product of a 1-dimensional stratified
space with $\real$ (the local orbits of the nonsingular vector field). 
Compactness and finiteness of the stratification imply that  
$\Sigma$ is everywhere locally homeomorphic to the product of 
a $K$-pronged radial tree with $\real$.

If $\Sigma$ is one-dimensional, then by transverse homogeneity and 
compactness of $\Sigma$, it is a compact one-manifold --- a circle.

If $\Sigma$ is two-dimensional, then every point of $\Sigma$ is 
locally homeomorphic to a $K$-pronged radial tree cross $\real$, 
where $K$ may vary but is always nonzero. The non-manifold
points of $\Sigma$ are precisely those points where $K\neq 2$. This 
set must be invariant under the flow, otherwise the uniqueness 
theorem for the vector field is violated. Hence, by compactness of 
the stratification, the non-manifold set is a finite link $L$.

Consider the space $\Sigma'$ obtained by removing from $\Sigma$ 
a small open tubular neighborhood of $L$. Since $L$ is an invariant
set for the flow, the vector field may be perturbed in such a way as to
leave $\Sigma'$ invariant. As all of the non-manifold points of 
$\Sigma$ have been removed, $\Sigma'$ is a 2-manifold with boundary.
The perturbed vector field on $\Sigma'$ is nonsingular; thus, 
$\Sigma'$ consists of annuli and M\"obius bands (plus perhaps tori
which do not encounter the singular link). 
\qed

We conclude the proof of Theorem~\ref{thm_Main} with the following
theorem.

\begin{theorem}
\label{thm_Unknots}
Any nonsingular vector field on $S^3$ having a $C^\omega$ integral of 
motion must possess a pair of unknotted closed orbits. 
\end{theorem}
\pf
We will prove this theorem for the slightly larger class of 
{\em stratified} integrals
which are not necessarily $C^\omega$ but whose critical sets 
are finite, [Whitney] stratified, and of positive codimension 
(see \S\ref{sec_Fomenko} for details and extensions).
Induct upon $\kappa$ the number of non-Bott connected components in the 
inverse images of the critical value set of the integral. If there are 
no such sets, then the system is Bott and the theorem follows from 
Corollary~\ref{cor_CMAN}. 

Let $c$ denote a (transversally) 
degenerate critical value of $P$ and $\Sigma$ a connected 
component of $P\inv(c)$. Denote by $N(\Sigma)$ the connected 
component of $P\inv([c-\eps,c+\eps])$ containing $\Sigma$. For 
$\eps$ sufficiently small $N$ is well-defined up to isotopy.  
The boundary components, $T_k$, of $N$ are all in 
the inverse image of regular values of the integral: as such, each 
$T_k$ is an embedded closed surface in $S^3$ supporting a nonsingular 
vector field -- a 2-torus. 

Each boundary torus $T_k$ bounds a 
solid torus in $S^3$ on at least one side \cite[p.107]{Rol77}. 
Denote by ${\mathcal S}$ the set of 
boundary components of $N$ which 
bound a solid torus containing $\Sigma$. 
Denote by ${\mathcal S}_0\subset{\mathcal S}$ the subset of
bounding neighborhood tori which are unknotted in $S^3$. 

{\bf Case 1:} ${\mathcal S}-{\mathcal S}_0\neq\emptyset$

Denote by $V$ the nontrivially-knotted solid torus containing $\Sigma$.
Redefine the integral on $V$ in a $C^\infty$ manner so 
that there is a single (Bott) critical set on the nontrivially knotted 
core of $V$, reducing $\kappa$. By the induction hypothesis, there 
is a pair of unknotted closed curves, neither of which can be 
the core of $V$. Note that although the new integral 
is not necessarily $C^\omega$, this stratified Bott integral is 
sufficient to apply Corollary~\ref{cor_CMAN} and the induction hypothesis.

{\bf Case 2:} ${\mathcal S}=\emptyset$

In this case, the non-Bott component $\Sigma$ must have a neighborhood
$N$ such $S^3-N$ consists of a disjoint collection of solid 
tori. We may then place a round-handle decomposition (RHD) on 
$S^3$ as follows.

By Lemma~\ref{lem_Branched}, one can decompose $\Sigma$ into a 
finite number of critical circles (1-strata) to which are attached 
annuli and M\"obius bands (2-strata) in a way which satisfies the 
Whitney condition.  Place an RHD on $N$ by thickening 
up each 1-stratum to a round 0-handle. The annular and M\"obius 
2-strata then thicken up in $N$ to round 1-handles of orientable 
and nonorientable type respectively. Since all of the boundary 
components of $N$ bound solid tori on the exterior of $N$, we can 
glue in round 2-handles, completing the RHD of $S^3$. 

According to the previously cited results of Asimov and Morgan, there is a 
nonsingular Morse-Smale flow on $S^3$ which realizes the indexed cores of 
this RHD as the periodic orbit link. Hence, by Theorem~\ref{thm_Wada}
there is a pair of unknots among the cores of this RHD. The index-2
cores are all nontrivially knotted by assumption; hence, the 
unknots have index zero or one. If zero, then these cores are 
the invariant 1-strata for the original flow. If an index-1 core is
unknotted, then there exists an invariant 2-stratum (annulus or 
M\"obius band) which is unknotted. If the core of an embedded 
annulus is unknotted, then both boundary components (invariant 
1-strata in the flow) are unknotted. 
In the case where the core of the invariant M\"obius 2-stratum is 
unknotted, we show that there exists an unknotted flowline as well
(see the proof of Lemma~\ref{lem_Mobius} below). 
Hence, both unknotted RHD cores are realized by isotopic invariant 
curves of the original flow.

{\bf Case 3:} ${\mathcal S}={\mathcal S}_0\neq\emptyset$

Construct a round-handle decomposition of $S^3$ as in 
Case 2 --- this is possible since all the exterior regions can be made
into round 2-handles. As before, there must exist a pair of 
unknotted cores to the RHD. Any which are of index zero or index one 
correspond to unknotted invariant curves in the original flow, by 
the arguments of Case 2. Assume that $V$ corresponds to an unknotted 
round 2-handle: a component of $S^3-N$. Replace the integral on 
$S^3-V$ (which is an unknotted solid torus as well) to have a 
single unknotted core critical set. Then by the induction hypothesis, 
there must have been an unknotted invariant curve within $V$. Hence, 
each unknotted round 2-handle corresponds to an unknotted invariant 
curve in the original flow.
\qed

\begin{lemma}
\label{lem_Mobius}
Any nonsingular vector field on a M\"obius band has a periodic orbit
isotopic to the core. 
\end{lemma}
\pf
The Poincar\'e-Bendixson Theorem holds for the Klein bottle (and thus for 
the M\"obius band which is a subset) by the theorem of Markley \cite{Mar69}.
The boundary curve of the M\"obius band is invariant, and either 
this curve has nontrivial holonomy or it does not. If the holonomy 
is nontrivial, then index theory and the Poincar\'e-Bendixson 
Theorem imply the existence of another closed orbit which is 
either twice-rounding (in which case is separates a smaller 
invariant M\"obius band --- repeat the analysis) 
or is once-rounding, in which case it 
is isotopic to the core. For the case of trivial holonomy, there is
a 1-parameter family of twice-rounding invariant curves, which 
either limits onto a closed curve with nontrivial holonomy, or 
else limits onto a once-rounding invariant core curve. 
\qed

We note that it is not necessarily the case that a stratified integral
on $S^3$ must have an unknotted curve of critical points (as is 
the case for a Bott integral). The centers of the M\"obius 
2-strata may be the only unknotted orbits in the flow: one may
construct an example in a manner reminiscent of a Seifert-fibred
structure on $S^3$ in which the critical sets of the integral are 
a pair of $(2,2n+1)$ torus knots whose (unknotted) cores are 
arranged in a Hopf link.

\section{Stratified integrable systems}
\label{sec_Fomenko}

The argument repeatedly stated in the Fomenko programme for 
restricting attention to Bott-integrable Hamiltonian systems is that 
this condition appears to be ubiquitous in physical
integrable systems (\ie, ones in which the integrals can be 
written out explicitly) \cite{Fom91,Fom88}. However, the examples cited
as evidence for this hypothesis are often real-analytic integrals. 
Hence, it would appear sensible to recast 
the Fomenko program of topological classification of integrable
Hamiltonian systems in the analytic case. 
Or, better still, one could allow for integrals 
which are less smooth yet satisfy the following more general 
conditions: one might call such integrals {\sc stratified}.
\begin{enumerate}
\item
	The critical values of the integral are isolated;
\item
	The inverse images of critical values are [Whitney] stratified 
	sets of codimension greater than zero.
\end{enumerate}
These assumptions allow for a controlled degeneracies in 
the integral, yet are by no means unnatural: codimension-one 
bifurcations of ``physical'' integrable Hamiltonian systems 
can and do exhibit such degeneracies. Since the topological results
of the Fomenko program can be obtained by using pre-existing 
RHD-theory \cite{CMAN}, and since the case of stratified 
non-Bott singularities also reduces to RHD's (see the proof 
of Theorem~\ref{thm_Unknots} above), there is seemingly no
reason to exclude stratified integrals. 

Many, if not all, of the key results of Fomenko's programme hold 
for this larger class of integrals. Unlike the Bott-Morse 
condition, it is often trivial to check whether the above 
criteria are met in the case of an explicit integral (as these
are almost always analytic functions).  We summarize a few  
results pertaining to the topology of flowlines which 
hold for stratified integrable systems. 
First, however,  we
recall that a {\sc graph-manifold} is a 3-manifold obtained 
by gluing together Seifert-fibred manifolds with boundary along
mutually incompressible tori \cite{Wa67a,Wa67b}. 

\begin{theorem}\label{thm_Stratified}
Given a nonsingular flow on a closed three-manifold $M$ possessing 
a stratified integral $P$, then $M$ is a graph-manifold. Furthermore,
if $M=S^3$, then the following statements hold:
\begin{enumerate}
\item
	There exists a pair of unknotted flowlines; 
\item
	Every closed orbit of the flow is a 
	knot which belongs to the family of zero-entropy knots
	described in Theorem~\ref{thm_Wada}.
\item
	The critical point set is {\sc nonsplittable}:
	there does not exist an embedded $S^2$ which separates 
	distinct components of $cp(P)$. 
\end{enumerate}
\end{theorem} 
\pf
The fact that $M$ is 
a graph-manifold follows trivially from the proof of Lemma~\ref{lem_Branched}
upon noting that a neighborhood of a two-dimensional degenerate set
has the structure of a round-handle decomposition. The results of 
Morgan \cite{Mor78} then imply that $M$ is a graph-manifold. Or, 
equivalently, one may perturb the integral on neighborhoods of 
singular sets to be Bott-Morse without changing the topology of 
the underlying manifold. This result was stated in \cite[p. 325]{Fom88}.

{\em Proof of Item (1):} 
Item 1 follows from the proof of Theorem~\ref{thm_Unknots} above.
\qed$_1$

{\em Proof of Item (2):} 
The following classification of knotted periodic orbits in stratified
integrable dynamics is an extension of the theorems of Casasayas 
et al. \cite{CMAN} and  Fomenko and Nguyen \cite{FN91} in the Bott
case. Choose any periodic orbit $\gamma$ whose $P$-value is not 
critical. Then, there exists a suitably small neighborhood of 
the non-Bott singular sets not containing $\gamma$. Perturb the 
integral on this small neighborhood to be a Bott integral. Although
this changes the vector field near the singular sets, it does not 
alter the knot type of $\gamma$, which must be a zero-entropy knot
by the aforementioned results.

If $\gamma$ is lying on a 2-stratum of the singular set 
then one may push $\gamma$ or a 2-cable of $\gamma$ off into a regular
torus $T$.  This regular torus is also a regular torus in
a Bott integral.  We may alter this integral to another Bott integral
for which the solid
torus that $T$ bounds contains a single critical level at the core.
Since the core of this torus is a zero entropy knot we know 
that $\gamma$ or a 2-cable of $\gamma$ is also a zero entropy
knot. Thus $\gamma$ is a zero entropy knot (see below).

The only case left to consider is a periodic orbit $\gamma$ lying on
the 1-stratum of the singular set. 
Some cable of $\gamma$ is a knot on the 2-stratum of the singular set and
we have merely to show, then, that if $\kappa$ is a zero-entropy knot
which is a cable of $\gamma$, then $\gamma$ is also zero-entropy.

Recall that zero-entropy knots are the closure of the unknot under
connected sum and cabling. If $\kappa$ is also a cable of a zero-entropy 
knot $\gamma'$, then we claim that $\gamma$ and $\gamma'$ are isotopic. 
Let $T$ and $T'$ denote 
the cabling tori for $\gamma$ and $\gamma'$ respectively. By transversality, 
$T\cap T'$ consists of disjoint circles having $\kappa$ as a component.
Any nullhomotopic circles can be inductively removed, leaving a 
finite collection of intersection curves isotopic to $\kappa$. These 
slice $T$ and $T'$ into pairs of annuli attached along their 
boundaries pairwise to form tori. One then uses the solid tori these 
bound to inductively cancel intersection curve pairs. Hence
the cores $\gamma$ and $\gamma'$ are isotopic.

In the other possibility, where $\kappa$ is the connected sum of two 
nontrivial knots, one has a contradiction upon showing that the 
nontrivial cable of a knot is always prime. A proof of this 
fact may be obtained by a similar geometric argument as the previous
step, or by an algebraic argument in \cite[p.93]{BZ85}. 
Hence, any periodic orbit is always a zero-entropy knot.
\qed$_2$

{\em Proof of Item (3):} 
Item (3) is seen to be true for Bott-integrable systems without critical 
tori by analyzing the operations of Theorem~\ref{thm_Wada} 
(see \cite{CMAN}). In the presence of critical tori, one can perturb the 
integral to have critical curves on the torus which renders the 
[now smaller] critical set unsplittable; hence the full critical set
was unsplittable as well.
In the stratified case, assume that $S$ is an embedded 2-sphere which 
separates the critical point set $cp(P)$. Then there exists a 
bound such that all sufficiently small smooth perturbations to the 
integral do not create critical points along $S$. Applying such a 
perturbation to a neighborhood of the inverse image of the 
critical values yields a Bott system with $S$ as a splitting 
sphere for the critical points set: contradiction.
\qed$_3$

This completes the proof of the Theorem. 
\qed

These results are noteworthy in that the existence of a 
single hyperbolic knot (\eg, a figure-eight knot) in a nonsingular 
vector field on $S^3$ implies the nonexistence
of an integral.


\end{document}